\documentclass{commat}

\usepackage{textcomp}

\title{%
    Cantor series expansions of rational numbers
   }

\author{%
  Symon Serbenyuk}

\affiliation{%
    45~Shchukina St.,  Vinnytsia, 21012, Ukraine.
    \email{simon6@ukr.net}
    }

\abstract{%
    This  survey  is devoted to necessary and sufficient conditions for a rational number to be representable by a Cantor series.  Necessary and sufficient conditions are formulated for the case of an arbitrary sequence $(q_k)$.  
   }



\VOLUME{31}
\YEAR{2023}
\NUMBER{1}
\firstpage{393}
\DOI{https://doi.org/10.46298/cm.10454}

\begin{paper}

\section{Introduction}

Let $Q\equiv (q_k)$ be a fixed sequence of positive integers, $q_k>1$,  $\Theta_k$ be a sequence of the sets $\Theta_k\equiv\{0,1,\dots ,q_k-1\}$, and $\varepsilon_k\in\Theta_k$.

The Cantor series expansion 
\begin{equation}
\label{eq: Cantor series}
\frac{\varepsilon_1}{q_1}+\frac{\varepsilon_2}{q_1q_2}+\dots +\frac{\varepsilon_k}{q_1q_2\dots q_k}+\dots
\end{equation}
of $x\in [0,1]$,   first studied by G. Cantor in \cite{Cantor1}, is a natural generalization of the b-ary expansion
$$
\frac{\alpha_1}{b}+\frac{\alpha_2}{b^2}+\dots+\frac{\alpha_n}{b^n}+\dots
$$
of numbers  from the closed interval $[0,1]$. Here $b$ is a fixed positive integer, $b>1$, and $\alpha_n\in\{0,1,\dots , b-1\}$.

By $x=\Delta^Q _{\varepsilon_1\varepsilon_2\ldots\varepsilon_k\ldots}$  denote a number $x\in [0,1]$ represented by series \eqref{eq: Cantor series}. This notation is called \emph{the representation of $x$ by Cantor series \eqref{eq: Cantor series}.}

We note that certain numbers from $[0,1]$ have two different representations by Cantor series \eqref{eq: Cantor series}, i.e., 
$$
\Delta^Q _{\varepsilon_1\varepsilon_2\ldots\varepsilon_{m-1}\varepsilon_m000\ldots}=\Delta^Q _{\varepsilon_1\varepsilon_2\ldots\varepsilon_{m-1}[\varepsilon_m-1][q_{m+1}-1][q_{m+2}-1]\ldots}=\sum^{m} _{i=1}{\frac{\varepsilon_i}{q_1q_2\dots q_i}}.
$$
Such numbers are called \emph{$Q$-rational}. The other numbers in $[0,1]$ are called \emph{$Q$-irrational}.

Cantor series expansions have been intensively studied from different points of view during the last century.
The metric, probability, and fractal theories of  number representations  by positive Cantor series were studied by a number of researchers. Also, functions  and fractal sets  defined in terms of Cantor series expansions were investigated. These problems were considered by the following researchers: P.~Erd\"os, J.~Galambos,  G.~Iommi, P.~Kirschenhofer, T.~Komatsu, V.~Laohakosol, B.~Li,   M.~Pa\v{s}t\'eka,  S.~Prugsapitak, J.~Rattanamoong, A.~R\'enyi, B.~Skorulski, R.~F.~Tichy, P. Tur\'an, Yi Wang, M.~S.~Waterman, H.~Wegmann, Liu Wen, Zhixiong Wen, Lifeng Xi, and other mathematicians.  

Such investigations  can be divided into two groups. The first is the investigation of  the fractional parts of  real numbers represented by Cantor series \eqref{eq: Cantor series}, and the other is the investigation of representations of non-negative integers represented by  positive Cantor series of the form
$$
n=\sum^{\infty} _{k=1}{\varepsilon_kq_1q_2\dots q_k},
$$
where  $\varepsilon_k\in \Theta_{k}$.

We give a brief description of these investigations.

A number of researches are devoted to studying various types of the normality of numbers represented by the Cantor series. In these papers, the notions of  Q-distribution normality, Q-normality, and Q-ratio normality,  are studied.  For example, in  the papers  \cite{Beros2014}, \cite{Mance_2012}, \cite{Mance_2014}, the notion of  Q-distribution normality is investigated. Indeed, one can note the following investigations: relations between various types of normality (e.g., see \cite{Beros2014}, \cite{Bill10}); the average value of the function of the sum of digits in the  Cantor series representation of a number  (see \cite{Kirschenhofer_Tichy84} and references in the last-mentioned article);  behaviour of the frequency of the most frequently used digit among the first digits in the representation of a number (e.g., see \cite{Erdos_Renyi1959}); necessary, sufficient, necessary and sufficient conditions for a number to be  a number having the property of certain type  normality   (see~\cite{Mance_2012}, \cite{Bill10}, \cite{Mance_2014});  the completeness of the Lebesgue measure, the density, topological properties, the Hausdorff measure of a set whose elements are numbers having the property of the normality  of a certain type (e.g., see~\cite{Mance_2012}, \cite{Mance_2014}); 
 the rationality and  irrationality of a number  which has  the property of the normality of a certain type (see \cite{Bill10}), etc.
Note that, in the papers \cite{Erdos_Renyi1959}, \cite{Erdos_Renyi1959(2)}, \cite{Renyi55}, \cite{Renyi56}, \cite{Renyi58}, \cite{Turan56},   P.~Erd\"os, A.~R\'enyi, and P. Tur\'an  introduced and studied the problem on 
normal numbers and other statistical properties of real numbers with respect to large classes of Cantor series expansions.
 Some investigations of Cantor series expansions were published by J.~Galambos in \cite{Galambos1976}, \cite{Galambos1976(2)}.

In some papers, certain generalizations of  real numbers representations  by the Cantor series are studied. For example,   properties of digits (sequences of digits) of the polyadic number $\alpha$ as functions (sequences of functions) of $ \alpha $ are studied in~\cite{Pasteka1996}; in  \cite{KLPR????}, the notion of a complex Cantor series is introduced, and  the $\mathbb Q$-algebraic and $\mathbb Q$-linearly independence of numbers represented by Cantor series are investigated; matrix expansions are studied in 
\cite{Waterman1975}; the papers \cite{Serbenyuk2016}, \cite{Serbenyuk2018} are devoted to certain generalizations of alternating Cantor series.

In certain papers, fractal properties of  representations of real numbers by positive Cantor series and fractal properties of certain type sets whose elements are represented by a positive Cantor series, are  studied (e.g., see~\cite{Iommi_Skorulski_2009}, \cite{Wang_Wen_Xi2009}, \cite{Mance_2015}, \cite{Bill10}). For example, in~\cite{Iommi_Skorulski_2009}, the Hausdorff-Besicovitch dimensions of sets whose elements are defined in terms of the frequencies of digits, are investigated. 
 The paper~\cite{Wegmann1968}  is devoted to studying  the conditions  under which the family of all possible rank cylinders $\Delta^Q _{c_1c_2\ldots c_n}$  is faithful for the Hausdorff-Besicovitch dimension calculation.
Sets whose elements have a restriction on using  digits in their own representations  are studied in \cite{Mance_2015}. In the last-mentioned article, the  formula for a calculation of the Hausdorff dimension of the following set is proved, and conditions for  the equality of the Hausdorff, packing, and box dimensions of this set, are discovered:
$$
R_I(Q)=\left\{x: x=\sum^{\infty} _{k=1}{\frac{\varepsilon_k}{q_1q_2\dots q_k}}, \varepsilon_k\in I_k\subseteq \Theta_k\right\}.
$$
Here the condition 
$$
\lim_{k\to\infty}{\frac{\log{q_k}}{\log{q_1q_2\dots q_k}}}=0
$$
 holds.

Also, we can note several investigations of functions.  The arguments or values of these functions are defined by positive~\cite{Cantor1} or alternating~\cite{Serbenyuk2017} Cantor series. In~\cite{Liu Wen2001}, properties of the following function were investigated:
$$
u=f(x)=\sum^{\infty} _{k=1}{\frac{u_k}{k(k+1)}},
$$ 
where $u_1=1$ and for $k=1,2,\dots$
$$
u_{k+1}=\begin{cases}
-\frac{u_k}{k},&\text{if $\varepsilon_{k+1}=0$ but $\varepsilon_k\ne 0$,}\\
&\text{or if $\varepsilon_{k+1}=q_{k+1}-1$ but $\varepsilon_k\ne q_k-1$;}\\
u_k,&\text{otherwise.}
\end{cases}
$$
Here $x$ is represented by series~\eqref{eq: Cantor series}. This function is well-defined and continuous. Also, $u=f(x)$ is  nowhere differentiable when $q_k\ge 3$ for all $k=1,2,\dots ,$ and the condition  
$$
\lim_{k\to\infty}{\frac{q_1q_2\dots q_k}{k!}}=\infty
$$
 holds. 
The last-mentioned function is a function with a complicated local structure. Certain examples of functions with a complicated local structure are described in~\cite{S. Serbenyuk preprint2}, \cite{Symon2015}, \cite{Symon2017}.
In the paper~\cite{Mance_2015}, the following function are studied:
$$
\psi_{P,Q}(x)=\sum^{\infty} _{k=1}{\frac{\min(E_k,q_k-1)}{q_1q_2\dots q_k}},
$$
where
$$
x=E_0+\sum^{\infty} _{k=1}{\frac{E_k}{p_1p_2\dots p_k}},~
\varepsilon_0+\sum^{\infty} _{k=1}{\frac{\varepsilon_k}{q_1q_2\dots q_k}}.
$$ 
 Here $E_0, \varepsilon_0\in\mathbb Z$,  $Q\equiv(q_k)$ and $P\equiv(p_k)$ are sequences of positive integers,  that greater than  $1$. Also, $E_k\ne p_k-1$ and $\varepsilon_k\ne q_k-1$ infinitely often, $E_k\in \{0,1,\dots, p_k-1\}$ and $\varepsilon_k\in \Theta_k$.

In the present article, the main attention is given to necessary and
sufficient conditions for $x$ (represented by  Cantor series with an arbitrary basic sequence $(q_k)$) to be rational. 

\begin{remark}
In the present article, we use  the following notations: $\mathbb N$, $\mathbb Z_0$, $\mathbb Z$, $\mathbb Q$, and $\mathbb I$.
Here by  $\mathbb N$  denote the set of all positive integers and by $\mathbb Z_0$ denotes the set $\mathbb N\cup \{0\}$, $\mathbb Z$ is the set of all integers, and $\mathbb Q$ is the set of all rational numbers, and $\mathbb I$ is the set of all irrational numbers.
\end{remark}

\section{Description of research of the main problem}

The problem of expansions of  rational/irrational numbers  in terms of generalizations of the b-ary numeral system is difficult. A version of this problem for Cantor series \eqref{eq: Cantor series} was introduced in the paper \cite{Cantor1} in 1869 and has been studied by a number of researchers.  For example,  G.~Cantor, P.~A.~Diananda, A.~Oppenheim, P.~Erd\"os, J.~Han\v{c}l, E.~G.~Straus, P.~Rucki, R.~Tijdeman, P.~Kuhapatanakul, V.~Laohakosol, D.~Marques, Pingzhi Yuan and other  scientists  studied this problem. 

In the monograph~\cite{Galambos1976}, Prof.~J\'anos~Galambos  called  the problem on representations of rational numbers by Cantor series \eqref{eq:  Cantor series} as  \emph{the fourth open problem}, and wrote the following: 

``Problem Four. Give a criterion of rationality for numbers given by a Cantor series. What one should seek here is a directly applicable criterion. A general sufficient condition for rationality would also be of interest, in which the quoted theorems of Diananda and Oppenheim (including the abstract criterion by condensations) can be guides or useful tools.

If in a Cantor series, negative and positive terms are permitted, somewhat less is known about the rationality or irrationality of the resulting sum. G. Lord (personal communication) tells me that the condensation method can be extended to this case as well, but still, the results are less complete than in the case of ordinary Cantor series."(\cite[p. 134]{Galambos1976}).

The paper~\cite{Serbenyuk21} is devoted to the last-mentioned discussion  and to expansions of rational numbers by sign-variable Cantor series. For fullness, one can note the following result of Diananda and Oppenheim noted by  J. Galambos.  
\begin{theorem}[\cite{Diananda_Oppenheim1955}]  
\label{th: Diananda Oppenheim}
A necessary and sufficient condition that $x$ given by~\eqref{eq: Cantor series} shall be rational is this: coprime integers $h, k$, $0\le h\le k$, an integer $N$ and a condensation shall exist such that
$$
A_i=\frac{h}{k}(B_i-1)
$$
for all $i\ge N$. 
\end{theorem}
Here
$$
x=X=A_0+\frac{A_1}{B_1}+\frac{A_2}{B_1B_2}+\dots+\frac{A_n}{B_1B_2\cdots B_n}+\dots,
$$
where $A_0=\varepsilon_0$ is the integer part of $x$, 
$$
B_1=q_1q_2\cdots q_{i_1}, B_2=q_{i_1+1}q_{i_2+1}\cdots q_{i_1+i_2},\dots ,
$$
 and $B_i\ge 2$, $0\le A_i \le B_i-1$, 
$$
\frac{\varepsilon_1}{q_1}+\frac{\varepsilon_2}{q_1q_2}+\dots+\frac{\varepsilon_{i_1}}{q_1q_2\cdots q_{i_1}}=\frac{A_1}{B_1}.
$$

We begin with a brief description of investigations of rational numbers represented by the Cantor series.  

Much research \cite{Tijdeman_Pingzhi2002}, \cite{Hancl97}, \cite{Hancl_Tijdeman2004}, \cite{Hancl2002}, \cite{Hancl_Tijdeman2004(2)} has been devoted to necessary or/and sufficient conditions for a rational number to be representable by   Cantor series \eqref{eq: Cantor series} such that sequences $(q_k)$ and $(\varepsilon_k)$ are sequences of integers. In some papers (see \cite{Hancl_Tijdeman2004}, \cite{Tijdeman_Pingzhi2002},\cite{Erdos_Straus1974}, \cite{Hancl2002}, \cite{Bill10}), the case of Cantor series for which  sequences $(q_k)$ and $(\varepsilon_k)$ are sequences of integers and the condition $\mathbb Z \ni q_k>1$ holds for all $k\in \mathbb N$, is investigated.  However, the main problem of the present article is studied for the case of series \eqref{eq: Cantor series} (e.g., see \cite{Cantor1}, \cite{Diananda_Oppenheim1955}, \cite{Kuhapatanakul_Laohakosol2001}, \cite{Oppenheim1954}) and still for the case of Cantor series of a special type (e.g., see \cite{Hancl_Tijdeman2005}, \cite{Hancl_Tijdeman2010}, \cite{Hancl_Tijdeman2009}, \cite{Hancl_Tijdeman2004}). For example, in the papers \cite{Diananda_Oppenheim1955}, \cite{Hancl_Tijdeman2004(2)}, \cite{Erdos_Straus1968}, \cite{KN2016}, Ahmes series are considered. The last series is   Cantor series  \eqref{eq: Cantor series} for which $\varepsilon_k=const=1$ holds for all $k\in \mathbb N$.

In the papers \cite{Diananda_Oppenheim1955}, \cite{Hancl97}, \cite{Hancl2002}, \cite{Hancl_Tijdeman2004}, \cite{Hancl_Tijdeman2004(2)}, \cite{Hancl_Tijdeman2005},\cite{Erdos_Straus1974}, \cite{Oppenheim1954}, \cite{Tijdeman_Pingzhi2002},  necessary and sufficient conditions for a rational (irrational) number to be representable by a  Cantor series are studied, and sufficient conditions are investigated in the papers \cite{Erdos_Straus1974}, \cite{Diananda_Oppenheim1955}, \cite{Hancl_Tijdeman2004}, \cite{Kuhapatanakul_Laohakosol2001}, \cite{Oppenheim1954}, \cite{Tijdeman_Pingzhi2002}. Although much research has been devoted to the problem of  representations of rational (irrational) numbers by   Cantor series for which sequences $(q_k)$ and $(\varepsilon_k)$ are sequences of special types
(see \cite{Cantor1},\cite{Erdos_Straus1974}, \cite{Hancl97}, \cite{Hancl2002}, \cite{Hancl_Tijdeman2004}, \cite{Kuhapatanakul_Laohakosol2001}, \cite{Oppenheim1954},\cite{Tijdeman_Pingzhi2002}),  little is known about necessary and sufficient conditions of the rationality (irrationality) for the case of an arbitrary sequence $(q_k)$ (see \cite{Diananda_Oppenheim1955}, \cite{Hancl_Tijdeman2004}, \cite{Serbenyuk: Cantor series}, \cite{S13}, \cite{Serbenyuk2017},  \cite{Rational numbers 2018},  \cite{Serbenyuk21}, \cite{Tijdeman_Pingzhi2002}). 

Finally,  several papers (see \cite{Hancl_Rucki2006}, \cite{KN2016}, \cite{Kuhapatanakul_Laohakosol2001}, \cite{Tijdeman_Pingzhi2002}) were devoted 
to investigations of conditions of the rationality or  irrationality of numbers represented by series of the form $\sum^{\infty} _{k=1}{\frac{a_k}{b_k}}$. Furthermore, in \cite{Kuhapatanakul_Laohakosol2001}, a  necessary and sufficient condition of the rationality of the sum  $\sum^{\infty} _{k=1}{\frac{a_k(-1)^{k+1}}{b_k}}$ is proved for the case of certain properties which are satisfied by  
 sequences $(a_k)$ and  $(b_k)$.

Let us consider our problem more in detail.

\section{Cantor's investigations, finite expansions, and conditions for finite expansions of rational numbers}

 Let us begin with a consideration of the results presented  in the first paper on this topic (i.e., \cite{Cantor1}). In \cite{Cantor1}, G. Cantor proved a fact that an arbitrary number $x\in [0,1)$  is a rational number if and only if $(\varepsilon_k)$ is ultimately periodic under the condition when a sequence $(q_k)$ is periodic. In addition, one can note the following theorem which  necessity  was given in~\cite{Cantor1} with the other formulation and with a more complicated proof for the case of positive Cantor series. 
 \begin{theorem}
\label{theorem2}
A rational number $\frac{p}{r}$ has a finite  expansion by a positive or sign-variable Cantor series  if and only if  there exists a number $n_0$ such that
 $$
q_1q_2\dots q_{n_0} \equiv 0\pmod{r}.
$$
\end{theorem}
The interest in the last theorem can be explained  (\cite{Serbenyuk: Cantor series}, \cite{S13}, \cite{Serbenyuk2017}, \cite{Serbenyuk21}) by the fact that      there   exist certain sequences $(q_k)$ such that all rational numbers represented by Cantor series (positive or sign-variable) have  finite expansions.         
For example, all rational numbers represented by the following representations have finite expansions. 
$$
x=\Delta^{-(2k)} _{\varepsilon_1\varepsilon_2\ldots \varepsilon_k...}\equiv\sum^{\infty} _{k=1}{\frac{(-1)^k\varepsilon_k}{2\cdot 4\cdot 8\cdot \ldots \cdot 2k}},~\text{where}~\varepsilon_k\in\{0,1,\dots,2k-1\};
$$
$$
x=\Delta^{(k+1)!} _{\varepsilon_1\varepsilon_2\ldots\varepsilon_k\ldots}\equiv\sum^{\infty} _{k=1}{\frac{\varepsilon_k}{2\cdot 3\cdot 4\cdot \ldots \cdot (k+1)}},~\text{where}~\varepsilon_k\in\{0,1,\dots,k\}.
$$

It is easy to see that there exist  sequences  $(q_k)$ and $(\varepsilon_k)$ such that a finite expansion is a necessary or/and sufficient condition of the rationality of any number represented by a Cantor series.  Several papers were devoted to such investigations. For example, see  \cite{Kuhapatanakul_Laohakosol2001}, \cite{Hancl2002}. Let us consider several related results.

In 2006, J. Sondow gave a geometric proof of the irrationality of the number~$e$~\cite{Sondow2006}. In \cite{Marques2009},  the following statement was proved by a generalization to  Sondow's construction. 
\begin{theorem}[{\cite{Marques2009}}]
Let $x=\sum^{\infty} _{k=1}{\frac{\varepsilon_n}{q_1q_2\cdots q_k}}$. Suppose that each prime divides
infinitely many of the $q_k$.  Then $x\in\mathbb I$ if and only if    both $0<\varepsilon_k<q_k-1$ hold infinitely often. 
\end{theorem}

For example, in~\cite{Hancl2002}, attention is given to conditions of finite expansions of rational numbers by positive and sign-variable Cantor expansions. That is, $\sum^{\infty} _{n=1}{\frac{\varepsilon_k}{q_1q_2\cdots q_k}}\in\mathbb Q$ if and only if  $\varepsilon_k=0$  for every sufficiently large positive integer $k$ under one of  the following two systems of conditions: 
\begin{itemize}
\item System 1 of conditions (the case of sign-variable series): suppose $(q_k)$ is  a sequence of positive integers greater than one, $(\varepsilon_k)$ is a sequence of integers such that the condition 
$$
\lim\inf_{k\to\infty}{\frac{|\varepsilon_k|+1}{q_k}}=0
$$
holds  and for every sufficiently large positive integer $k$
$$
|\varepsilon_{k+1}|\le \frac{1}{2}\max{(|\varepsilon_k|,1)}q_{k+1}.
$$
\item System 2 of conditions (the case of positive series): suppose $(q_k)$ is a sequence of positive integers greater than one   and $K\in(0,1)$, $(\varepsilon_k)$ is a sequence of non-negative integers such that the condition
$$
\lim\inf_{n\to\infty}{\frac{\varepsilon_k+1}{q_k}}=0
$$
holds   and for every sufficiently large positive integer $k$
$$
\varepsilon_{k+1}\le K\max{(\varepsilon_k,1)}q_{k+1}.
$$
\end{itemize}

\section{The shift operator and related investigations}

We must note that the notion of the shift operator plays an important role in investigations of expansions of rational numbers defined by the Cantor series (positive, alternating, or sign-variable).

We begin with definitions. Let $   N_B$ be a fixed subset of positive integers, 
$$
\rho_k=\begin{cases}
1&\text{if $k\in  N_B$}\\
2&\text{if  $k\notin  N_B$,}
\end{cases}
$$
and $Q\equiv (q_k)$ be a fixed sequence of positive integers such that $q_k>1$ for all $n\in\mathbb   N$.  Then we get the following representation of real numbers
\begin{equation}
\label{eq: sign-variable series}
x=\Delta^{(\pm Q,  N_B)} _{\varepsilon_1\varepsilon_2\ldots\varepsilon_k\ldots}\equiv\frac{(-1)^{\rho_1}\varepsilon_1}{q_1}+\frac{(-1)^{\rho_2}\varepsilon_2}{q_1q_2}+\dots +\frac{(-1)^{\rho_k}\varepsilon_k}{q_1q_2\dots q_k}+\dots, 
\end{equation}
where $\varepsilon_k\in\{0,1,\dots , q_k-1\}$. 

The last representation is called \emph{the representation of a number $x$ by a sign-variable Cantor series} or \emph{the quasi-nega-Q-representation}.  It is easy to see that we get a positive Cantor series whenever $N_B=\emptyset$.

 Define \emph{the shift operator $\sigma$ of expansion \eqref{eq: sign-variable series}} by the rule
$$
\sigma(x)=\sigma\left(\Delta^{(\pm Q, N_B)} _{\varepsilon_1\varepsilon_2\ldots\varepsilon_k\ldots}\right)=\sum^{\infty} _{k=2}{\frac{(-1)^{\rho_k}\varepsilon_k}{q_2q_3\dots q_k}}=q_1\Delta^{(\pm Q, N_B)} _{0\varepsilon_2\ldots\varepsilon_k\ldots}.
$$
Clearly,
\begin{equation}
\label{eq:  series 2}
\begin{split}
\sigma^n(x) &=\sigma^n\left(\Delta^{(\pm Q, N_B)} _{\varepsilon_1\varepsilon_2\ldots\varepsilon_k\ldots}\right)\\
& =\sum^{\infty} _{k=n+1}{\frac{(-1)^{\rho_k}\varepsilon_k}{q_{n+1}q_{n+2}\dots q_k}}=q_1\dots q_n\Delta^{(\pm Q, N_B)} _{\underbrace{0\ldots 0}_{n}\varepsilon_{n+1}\varepsilon_{n+2}\ldots}.
\end{split}
\end{equation}

The following theorem is the most general statement on the representation of rational numbers for any sequences $(q_k)$,  $(\varepsilon_k)$, and an arbitrary set $N_B$.
\begin{theorem}[\cite{Serbenyuk: Cantor series}, \cite{S13}, \cite{Serbenyuk2017}, \cite{Rational numbers 2018}]
\label{th: the main theorem}
A number $x$ represented by series \eqref{eq: sign-variable series} is  rational for the case of any $N_B \subseteq~\mathbb N$ if and only if  there exist numbers $n\in\mathbb Z_0$ and $m\in\mathbb N$ such that $\sigma^n(x)=\sigma^{n+m}(x)$.
\end{theorem}
The last theorem can be formulated by the following way. 
\begin{theorem}[\cite{Serbenyuk: Cantor series}, \cite{S13}, \cite{Rational numbers 2018}]
\label{th: the main theorem 2}
A number $x=\Delta^{(\pm Q, N_B)} _{\varepsilon_1\varepsilon_2\ldots \varepsilon_k\ldots }$  is  rational if and only if there exist numbers $n\in\mathbb Z_0$ and $m\in\mathbb N$ such that
$$
\Delta^{(\pm Q, N_B)} _{\underbrace{0\ldots 0}_{n}\varepsilon_{n+1}\varepsilon_{n+2}\ldots }=q_{n+1}\dots q_{n+m}\Delta^{(\pm Q, N_B)} _{\underbrace{0\ldots 0}_{n+m}\varepsilon_{n+m+1}\varepsilon_{n+m+2}\ldots }.
$$
\end{theorem}

Let us recall several auxiliary statements which are true for positive Cantor series but do not hold for the general case of  sign-variable Cantor series (i.e., for certain sets $N_B$). 
\begin{lemma}[\cite{Serbenyuk: Cantor series}, \cite{S13}]
Let $n_0$ be a fixed positive integer. Then the condition $\sigma^n(x)=const$ holds for all $n\ge n_0$ if and only if  $\frac{\varepsilon_n}{q_n-1}=const$ for all $n>n_0$.
\end{lemma}
\begin{lemma}[\cite{Serbenyuk: Cantor series}, \cite{S13}]
Suppose we have $q=\min_{n\in\mathbb N} {q_n}$ and  fixed $\varepsilon\in\{0,1,\dots ,q-1\}$. Then the condition $\sigma^n(x)=x=\frac{\varepsilon}{q-1}$ holds if and only if the condition $\frac{q_n-1}{q-1}\varepsilon=\varepsilon_n\in\mathbb Z_0$ holds for all $n\in\mathbb N$.
\end{lemma}

Let us consider cases when the condition  $\frac{\varepsilon_k}{q_k-1}=const$ (the last equality holds for all $k$ greater than some fixed $k_0$)  is a  necessary and/or  sufficient condition for a rational number to be representable by a  positive Cantor series. For more information, see \cite{Diananda_Oppenheim1955}, \cite{Hancl_Tijdeman2004}, \cite{Tijdeman_Pingzhi2002}.  

In  \cite{Hancl_Tijdeman2004},  J.~Han\v{c}l and  R.~Tijdeman  formulated certain conditions of the irrationality of a number represented by  Cantor series \eqref{eq: Cantor series} when sequences $(q_k)$ and $(\varepsilon_k)$  are sequences of positive integers and 
$q_k>1$ for all $k\in \mathbb N$. Applications of the shift operator to representations of rational numbers by such  series are considered.  This article is partially devoted to 
 conditions under which the condition $\frac{\varepsilon_k}{q_k-1}=const$ is a necessary and sufficient condition of the rationality of  numbers represented by such expansions. In particular, the following cases are considered: \[
    \lim\inf_{k\to\infty}\left(\frac{\varepsilon_{k+1}}{q_{k+1}}-\frac{\varepsilon_k}{q_k}\right)=0, \quad
    \varepsilon_k=o(q_{k-1}q_k), \quad
    \varepsilon_{k+1}-\varepsilon_k=o(q_{k-1}q_k).
\]
Also, in \cite{Hancl_Tijdeman2004}, the authors noted that sum  \eqref{eq: Cantor series} is equal to a rational number if  
$\frac{\varepsilon_k}{q_k-1}=const$ holds for all  $k$ greater than some number $n_0$. Let us recall some results. 
 \begin{lemma}[\cite{Hancl_Tijdeman2004}]
If $S=\sum^{\infty} _{k=1}{\frac{\varepsilon_k}{q_1q_2\ldots q_k}}=\frac{r}{p}$ holds for a certain   $r\in \mathbb Z$ and  $p\in\mathbb N$, then $pS_N\in\mathbb Z$ for all  $N\in\mathbb N$.
\end{lemma}
Here $S=\sigma^0 (x)$ and $S_N=\sigma^{N-1} (x)$. That is, 
$$
S_N=\sum^{\infty} _{k=N}{\frac{\varepsilon_n}{q_N\cdots q_k}}.
$$
\begin{proposition}[\cite{Hancl_Tijdeman2004}]
If $(S_k)$ is bounded from below and for every   $\varepsilon>0$ we have  
$$S_{k+1}-S_k<\varepsilon
$$ for  $k\ge k_0(\varepsilon)$, then $S=\sum^{\infty} _{k=1}{\frac{\varepsilon_k}{q_1q_2\cdots q_k}}\in \mathbb Q$  if and only if  $\frac{\varepsilon_k}{q_k-1}=const$ for $N>N_0$.
\end{proposition}
\begin{corollary}[\cite{Hancl_Tijdeman2004}]
If $(\varepsilon_k)$ is a sequence of positive integers such that  $\varepsilon_{k+1}-\varepsilon_k=o(k)$, then  $\sum^{\infty} _{k=1}{\frac{\varepsilon_k}{k!}}\in\mathbb Q$ if and only if $\frac{\varepsilon_k}{k-1}=const$ for $k$ greater than some $k_1$.
\end{corollary}
\begin{theorem}[{{\cite{Hancl_Tijdeman2004}}}]
Let $(q_k)$ be a sequence of positive integers which is monotonic and satisfies $\varepsilon_k=o(q^2 _k)$. Then $\sum^{\infty} _{k=1}{\frac{\varepsilon_k}{q_1q_2\cdots q_k}}\in \mathbb Q$  if and only if  $\frac{\varepsilon_k}{q_k-1}=const$ for $k\ge k_0$.
\end{theorem}
\begin{theorem}[{{\cite{Hancl_Tijdeman2004}}}]
Let  $(q_k)$ and $(\varepsilon_k)$ be sequences of integers such that  $q_k>1$ for all $k\in \mathbb N$. If $\left(\frac{\varepsilon_k}{q_k}\right)$ is  bounded from below,   $\lim_{k\to \infty}{\frac{\varepsilon_k}{q_{k-1}q_k}}=0$, and 
for each $\varepsilon>0$ there exists $k_0(\varepsilon)$ such that the condition $\frac{\varepsilon_{k+1}}{q_{k+1}}<\frac{\varepsilon_k}{q_k}+\varepsilon$  holds for  $k>k_0(\varepsilon)$, then   $\sum^{\infty} _{k=1}{\frac{\varepsilon_k}{q_1q_2\cdots q_k}}\in \mathbb Q$  if and only if  $\frac{\varepsilon_k}{q_k-1}=const$ for $k\ge N_0$.
\end{theorem}
\begin{theorem}[{{\cite{Hancl_Tijdeman2004}}}]
Let $(q_k)$ be a monotonic sequence of positive integers satisfying   $\lim_{k\to\infty}{\frac{q_k}{\log k}}=\infty$. Then $\sum^{\infty} _{k=1}{\frac{\varepsilon_k}{q_1q_2\cdots q_k}}\in \mathbb Q$  if and only if  $\frac{\varepsilon_k}{q_k-1}=const$ for $k\ge k_0$.
\end{theorem}
\begin{theorem}[{{\cite{Hancl_Tijdeman2004}}}]
Let  $(q_k)$ be an unbounded monotonic sequence of positive integers. Then $\sum^{\infty} _{k=1}{\frac{k}{q_1q_2\cdots q_k}}\in \mathbb Q$  if and only if $\frac{k}{q_k-1}=const$ for $k\ge k_0$.
\end{theorem}

Results obtained in \cite{Hancl_Tijdeman2004} were generalized  and  corrected by  Robert Tijdeman and Pingzhi Yuan in paper 
\cite{Tijdeman_Pingzhi2002}. In particular, results are  generalized for the cases when $\varepsilon_k=k$ and $q_k\to\infty$,   $q_k=k$ and $\varepsilon_{k+1}-\varepsilon_k=O(k)$.  In the last-mentioned article, it is shown that, in order that the condition $\frac{\varepsilon_k}{q_k-1}=const$ for all $k\ge k_0$ is a necessary and sufficient condition of the rationality,  one can  neglect the condition $\varepsilon_k=o(q^2 _k)$ in the system of conditions: $\varepsilon_k=o(q^2 _k)$ , $\varepsilon_k \ge 0$, $\varepsilon_{k+1}-\varepsilon_k<\varepsilon q_k$ for $k\ge  k_1(\varepsilon)$. We note the following statements.

\begin{theorem}[{{\cite{Tijdeman_Pingzhi2002}}}]
Let $(q_n)$ be a monotonic integer sequence with $q_n>1$ for all $n$ and  $(\varepsilon_n)$ be an integer sequence such that
$\varepsilon_{n+1}-\varepsilon_n=o(q_{n+1})$. Then $\sum^{\infty} _{n=1}{\frac{\varepsilon_n}{q_1q_2\cdots q_n}}\in \mathbb Q$ if and only if $\frac{\varepsilon_n}{q_n-1}=const$  for all $n$ greater than some $ n_0$.
\end{theorem}

\begin{theorem}[{{\cite{Tijdeman_Pingzhi2002}}}]
Let $(q_k)$ be a monotonic sequence of positive integers, $q_k>1$. Let $(\varepsilon_k)$ be a sequence of positive integers  satisfying 
$$
\lim\sup_{k\to\infty}{\frac{\varepsilon_{k+1}-\varepsilon_k}{q_k}}\le 0.
$$
Then  $\sum^{\infty} _{k=1}{\frac{\varepsilon_k}{q_1q_2\ldots q_k}}\in \mathbb Q$ if and only if  $\frac{\varepsilon_k}{q_k-1}=const$ for all $k$ greater than some  $k_0$.
\end{theorem} 
In addition, the following sufficient condition of the irrationality is proved. 
\begin{theorem}[{{\cite{Tijdeman_Pingzhi2002}}}]
Let $q_k>1$ be such that $\varepsilon_k=O(q_k)$ for all $k$ and
 $\lim_{k\to\infty}{\frac{\varepsilon_k}{q_k}}=\alpha \in\mathbb I$.  Then $\sum^{\infty} _{k=1}{\frac{\varepsilon_k}{q_1q_2\cdots q_k}}\in\mathbb I$.
\end{theorem}
The last statement with the condition  $0\le \varepsilon_k<q_k$ without $\varepsilon_k=O(q_k)$ was proved in \cite{Oppenheim1954}.

Finally, in  \cite{Tijdeman_Pingzhi2002}, the following denotations are used in proofs:
$$
S=\sum^{\infty} _{k=1}{\frac{\varepsilon^{*} _{k}}{q^{*} _1q^{*} _2\cdots q^{*} _k}},~ S_{n_k}=\sum^{k} _{j=1}{\frac{\varepsilon^{*} _j}{q^{*} _1q^{*} _2\cdots q^{*} _j}}, ~ R_{n_k}=\sum^{\infty} _{j=k+1}{\frac{\varepsilon^{*} _j}{q^{*} _{k+1}q^{*} _{k+2}\cdots q^{*} _j}}.
$$
Here  $(n_k)$ is a subsequence of positive integers, $n_0=1$,
$$
\varepsilon^{*} _k=\varepsilon_{n_k-1}+\varepsilon_{n_k-2}q_{n_k-1}+\dots+\varepsilon_{n_{k-1}}q_{n_k-1}q_{n_k-2}\cdots q_{n_{k-1}+1},
$$
and $q^{*} _k=q_{n_k-1}q_{n_k-2}\cdots q_{n_{k-1}}$,  $k=1,2,3,\ldots$.  
For series \eqref{eq: Cantor series}, where $(q_n)$ and $(\varepsilon_n)$ are sequences of integers such that $q_n>0$ for all $n\in\mathbb N$ and series \eqref{eq: Cantor series} converges, the following statements are true.  
\begin{lemma}[{{\cite{{Tijdeman_Pingzhi2002}}}}]
Using the notation above, if there exists a subsequence $(n_k)$ of  positive integers such that 
 $R_{n_k}=R_{n_{k+1}}$  for  $k=1,2,\dots$, then  $S\in\mathbb Q$.
\end{lemma}
\begin{proposition}[{\cite{Tijdeman_Pingzhi2002}}]
If  $(R_n)$ is bounded from below and there exists a subsequence $(n_k)$ of positive integers with $R_{n_{k+1}}-R_{n_k}<\varepsilon$ for  $k\ge k_0(\varepsilon)$, then $S\in\mathbb Q$ if and only if $R_{n_k}=R_{n_{k+1}}$ for all large
$k$.
\end{proposition}

In \cite{Oppenheim1954}, A. Oppenheim studied sufficient conditions of the irrationality of numbers represented by Cantor series \eqref{eq: Cantor series} and, also, alternating series  \eqref{eq: Cantor series} such that $|\varepsilon_i|<q_i-1$ for $i=1,2,3,\dots$, and $\varepsilon_m\varepsilon_n<0$ for some $m>i$ and $n>i$ when $i$ is any fixed integer. Also, in \cite{Oppenheim1954}, the main results obtained by using some results from  \cite{Cantor1} and  sums of the form
$$
x_{i_k}=\frac{\varepsilon_{i_k}}{q_{i_k}}+\frac{\varepsilon_{i_k+1}}{q_{i_k}q_{i_k+1}}+\frac{\varepsilon_{i_k+2}}{q_{i_k}q_{i_k+1}q_{i_k+2}}+\dots ,
$$
 where $(i_k)$ is some subsequence of positive integers, and by  investigation of the limit of $c_{i_k}=\frac{\varepsilon_{i_k}}{q_{i_k}}$ as $k\to\infty$. That is, here $x_{i_k}=\sigma^{i_k-1}(x)$.
 \begin{lemma}[\cite{Oppenheim1954}]
A necessary and sufficient condition that $x$ given by  convergent series  \eqref{eq: Cantor series}, where $q_k$ and  $\varepsilon_k$ are integers, shall be irrational is that for every integer $p\in\mathbb N$ we can find an integer  $r\in\mathbb Z$ 
and a subsequence $(i_k)$ such that
$$
\frac{r}{p}<x_{i_k}<\frac{r+1}{p},~k=1,2,3,\ldots .
$$ 
\end{lemma}

 Finally, in this section, we note necessary and sufficient conditions for a rational number to be representable by  certain types of  Cantor series which were  investigated by  P.~Erd\"os and E.~G.~Straus  in \cite{Erdos_Straus1974}.
\begin{theorem}[{{\cite{Erdos_Straus1974}}}]
Let  $(\varepsilon_n)$ be a sequence of integers and $(q_k)$ be a sequence of positive integers with $q_k>1$ for all large $k$ and 
$$
\lim_{k\to\infty}{\frac{|\varepsilon_k|}{q_{k-1}q_k}}=0.
$$
Then $\sum^{\infty} _{k=1}{\frac{\varepsilon_k}{q_1q_2\cdots q_k}}\in\mathbb Q$ if and only if there exist a positive integer
$B$  and a sequence of integers $(c_k)$  such that for all large $k$ we have
$$
B\varepsilon_k=c_kq_k-c_{k+1}, ~|c_{k+1}|<\frac{q_k}{2}.
$$
\end{theorem}
\begin{theorem}[{{\cite{Erdos_Straus1974}}}]
Let $p_k$ be the $k$th prime and let $(q_k)$  be a monotonic sequence of positive integers satisfying 
$$
\lim_{k\to\infty}{\frac{p_k}{q^2 _k}}=0, ~\lim\inf_{k\to\infty}{\frac{q_k}{p_k}}=0.
$$
Then $\sum^{\infty} _{k=1}{\frac{p_k}{q_1q_2\cdots q_k}}\in\mathbb I$.
\end{theorem}

\section{Certain approaches to investigations of expansions of rational numbers}

We considered mainly the shift operator.  Now one can consider  another approaches to investigations of expansions of rational numbers but some of them are related with the shift operator. 

 In \cite{{Bill10}},  the probabilistic  approach is used and the attention is given to Cantor series~\eqref{eq: Cantor series} for which  $\varepsilon_k\ne q_k-1$ infinitely often. Irrationality of numbers  having a property of a certain type of the normality is investigated. 
\begin{definition}[{\cite[p.45]{{Bill10}}}]
A number $x\in [0,1)$ is called  \emph{Q-distribution normal} if the sequence
$$
X=(x \pmod{1}, q_1x \pmod{1}, q_1q_2x \pmod{1}, q_1q_2\cdots q_kx \pmod{1}, \dots )
$$
 is uniformly distributed in $[0,1)$.
\end{definition}
\begin{theorem}[{\cite[p. 264]{{Bill10}}}]
A number $x\in[0,1)$ is irrational if and only if there exists a basic sequence  $Q=(q_k)$ such that $x$ is Q-distribution normal.
\end{theorem}

In the paper \cite{Kumar2019}, the subspace theorem is used for proving conditions for a transcendental number to be representable by positive Cantor series. Such conditions were formulated in terms of  blocks of digits $\varepsilon_k$  and in terms of tuples of digits for expansion~\eqref{eq: Cantor series}.

Finally, one approach  based on the notion of cylinders of Cantor expansions gives an opportunity to  model rational numbers. However, we have necessary and sufficient conditions  for the case of the positive Cantor series and  a necessary condition for the sign-variable Cantor series. Let us consider the following two theorems. 
\begin{theorem}[\cite{Rational numbers 2018}]
\label{th: modeling}
A number $x=\Delta^Q _{\varepsilon_1\varepsilon_2\ldots\varepsilon_n\ldots} \in (0,1)$ represented by series \eqref{eq: Cantor series}  is a rational number $\frac{p}{r}$, where $p,r\in\mathbb N, (p,r)=1$, and $p<r$, if and only if the condition 
$$
\varepsilon_n=\left[\frac{q_n(\Delta_{n-1}-r\varepsilon_{n-1})}{r}\right]
$$
 holds for all $1<n\in\mathbb N$, where $\Delta_1=pq_1$, $\varepsilon_1=\left[\frac{\Delta_1}{r}\right]$, and 
 $[a]$ is the integer part of $a$.
\end{theorem}
For fullness, we give some examples of rational numbers from~\cite{Rational numbers 2018}. Really, suppose
$$
x=\Delta^{(2n+1)} _{\varepsilon_1\varepsilon_2\ldots\varepsilon_n\ldots}
 =\sum^{\infty} _{n=1}{\frac{\varepsilon_n}{3\cdot 5\cdot 7 \cdots  (2n+1)}}  . 
$$
Then 
$$
 \frac{1}{4}=\Delta^{(2n+1)} _{035229[11]4\ldots},
 \qquad
 \frac{3}{8}=\Delta^{(2n+1)} _{104341967\ldots}.
$$
\begin{theorem}[\cite{Serbenyuk21}]
If $x=\Delta^{(\pm Q, N_B)} _{\varepsilon_1\varepsilon_2\ldots\varepsilon_k\ldots} =\frac{p}{r}$, where $p\in\mathbb Z, r\in\mathbb N, (|p|,|r|)=1$, and $|p|<r$, then the condition 
$$
\varepsilon_k=\left|\left[\frac{q_k(\Delta^{(g)} _{k-1}-(-1)^{\rho_{k-1}}r\varepsilon_{k-1})}{r}+s_n\right]\right|
$$
 holds for all $1<k\in\mathbb N$. Here $\Delta^{(g)} _1=pq_1$,  $\varepsilon_1=\left|\left[\frac{\Delta^{(g)} _1}{r}+s_1\right]\right|$, and
 $[a]$ is the integer part of $a$. Also, 
$$
s_1=\sum_{1<k \in N_B}{\frac{q_k-1}{q_2q_3\cdots q_k}},
\qquad
s_k=\begin{cases}
q_ks_{k-1}&\text{whenever $k\notin N_B$}\\
q_ks_{k-1}-(q_k-1)&\text{whenever $k\in N_B$.}
\end{cases}
$$
\end{theorem}

One can note that the two last statements are related to the shift operator. Really, in the case  of positive Cantor series \cite{Rational numbers 2018}, we have $\sigma^n(x)=\{\frac{\Delta_n}{r}\}$ and $\varepsilon_n=[\frac{\Delta_n}{r}]$. In the general  case of sign-variable Cantor series (i.e., there is no  number $k_0$ such that any $k\in N_B$ or any $k\notin N_B$ for all $k>k_0$), we obtain~\cite{Serbenyuk21} the following:
$$
\left\{\frac{\Delta^{(g)} _n}{r}\right\}=\begin{cases}
\sigma^n(x)&\text{whenever $\sigma^n(x)\ge 0$}\\
1- \sigma^n(x)&\text{whenever $\sigma^n(x)<0$,}
\end{cases}
$$
where  $\{a\}$ is the fractional part of $a$ (i.e., $a=[a]+\{a\}$).

In this survey, we have demonstrated the main  conditions for a rational number to be representable by   positive, alternating,  and 
sign-variable Cantor series. Connections among some of them are described. An important role of the notion of the shift operator in investigations in this topic, is noted.


\EditInfo{May 02, 2021}{February 21, 2022}{Attila Bérczes}

\end{paper}